\documentclass{amsart}
\usepackage{amsmath,amssymb,amsthm}
\usepackage[pdftex]{graphicx}
\usepackage{graphicx}
\usepackage{longtable}

\newtheorem{defi}{Definition}[section]
\newtheorem{theo}{Theorem}[section]
\newtheorem{lemm}[theo]{Lemma}
\newtheorem{prop}[theo]{Proposition}
\newtheorem{rem}[theo]{Remark}

\makeatletter
    
    \@addtoreset{equation}{section}
\makeatother
\allowdisplaybreaks[4]

\newcommand{\norm}[2]{\|#1\|_\mathfrak{s}^{#2}}
\newcommand{\supp}{\mathop{\text{\rm supp}}}
\newcommand{\id}{\text{\rm id}}
\newcommand{\ind}{\mathbf{1}}

\newcommand{\ig}[2]{\includegraphics[trim=0 #1 0 0]{#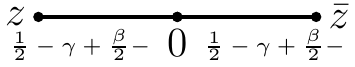}}

\newcommand{\f}[2]{\frac{#1}{#2}}
\newcommand{\cal}[1]{\mathcal{#1}}
\newcommand{\spa}[1]{\langle#1\rangle}

\def\${|\!|\!|}

\def\a{\alpha}
\def\b{\beta}
\def\d{\partial}
\def\e{\epsilon}
\def\g{\gamma}
\def\k{\kappa}
\def\s{\mathfrak{s}}
\def\t{\theta}

\def\E{\mathbb{E}}
\def\K{\mathfrak{K}}
\def\R{\mathbb{R}}
\def\T{\mathbb{T}}
\def\Z{\mathbb{Z}}

\def\vp{\varphi}

\newsavebox{\boxV}
\sbox{\boxV}{\includegraphics[trim=0 3 0 0]{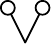}}
\def\V{\usebox{\boxV}}

\newsavebox{\boxmV}
\sbox{\boxmV}{\includegraphics[width=3mm, trim=0 3 0 0]{symbol_V.pdf}}
\def\mV{\usebox{\boxmV}}

\newsavebox{\boxW}
\sbox{\boxW}{\includegraphics[trim=0 3 0 0]{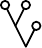}}
\def\W{\usebox{\boxW}}

\newsavebox{\boxI}
\sbox{\boxI}{\includegraphics[trim=0 2 0 0]{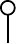}}
\def\I{\usebox{\boxI}}

\newsavebox{\boxB}
\sbox{\boxB}{\includegraphics[trim=0 5 0 0]{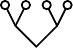}}
\def\B{\usebox{\boxB}}

\newsavebox{\boxmB}
\sbox{\boxmB}{\includegraphics[width=3mm, trim=0 5 0 0]{symbol_B.pdf}}
\def\mB{\usebox{\boxmB}}

\newsavebox{\boxWV}
\sbox{\boxWV}{\includegraphics[trim=0 7 0 0]{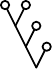}}
\def\WV{\usebox{\boxWV}}

\newsavebox{\boxmWV}
\sbox{\boxmWV}{\includegraphics[width=3mm, trim=0 7 0 0]{symbol_WV.pdf}}
\def\mWV{\usebox{\boxmWV}}

\newsavebox{\boxY}
\sbox{\boxY}{\includegraphics[trim=0 6 0 0]{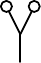}}
\def\Y{\usebox{\boxY}}

\newsavebox{\boxC}
\sbox{\boxC}{\includegraphics[trim=0 3 0 0]{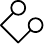}}
\def\C{\usebox{\boxC}}

\newsavebox{\boxmC}
\sbox{\boxmC}{\includegraphics[width=3mm, trim=0 3 0 0]{symbol_C.pdf}}

\newsavebox{\boxVB}
\sbox{\boxVB}{\includegraphics[trim=0 7 0 0]{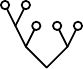}}
\def\VB{\usebox{\boxVB}}

\newsavebox{\boxBV}
\sbox{\boxBV}{\includegraphics[trim=0 7 0 0]{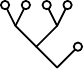}}
\def\BV{\usebox{\boxBV}}

\newsavebox{\boxWW}
\sbox{\boxWW}{\includegraphics[trim=0 10 0 0]{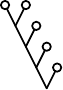}}
\def\WW{\usebox{\boxWW}}

\newsavebox{\boxF}
\sbox{\boxF}{\includegraphics[trim=0 5 0 0]{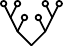}}
\def\F{\usebox{\boxF}}

\newsavebox{\boxmF}
\sbox{\boxmF}{\includegraphics[width=4mm, trim=0 5 0 0]{symbol_F.pdf}}
\def\mF{\usebox{\boxmF}}

\newsavebox{\boxBW}
\sbox{\boxBW}{\includegraphics[trim=0 5 0 0]{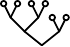}}
\def\BW{\usebox{\boxBW}}

\newsavebox{\boxmBW}
\sbox{\boxmBW}{\includegraphics[width=4mm, trim=0 5 0 0]{symbol_BW.pdf}}
\def\mBW{\usebox{\boxmBW}}

\newsavebox{\boxWB}
\sbox{\boxWB}{\includegraphics[trim=0 7 0 0]{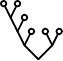}}
\def\WB{\usebox{\boxWB}}

\newsavebox{\boxmWB}
\sbox{\boxmWB}{\includegraphics[width=4mm, trim=0 5 0 0]{symbol_WB.pdf}}
\def\mWB{\usebox{\boxmWB}}

\newsavebox{\boxVBV}
\sbox{\boxVBV}{\includegraphics[trim=0 7 0 0]{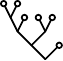}}
\def\VBV{\usebox{\boxVBV}}

\newsavebox{\boxmVBV}
\sbox{\boxmVBV}{\includegraphics[width=4mm, trim=0 5 0 0]{symbol_VBV.pdf}}
\def\mVBV{\usebox{\boxmVBV}}

\newsavebox{\boxBVV}
\sbox{\boxBVV}{\includegraphics[trim=0 7 0 0]{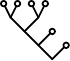}}
\def\BVV{\usebox{\boxBVV}}

\newsavebox{\boxmBVV}
\sbox{\boxmBVV}{\includegraphics[width=4mm, trim=0 5 0 0]{symbol_BVV.pdf}}
\def\mBVV{\usebox{\boxmBVV}}

\newsavebox{\boxWWV}
\sbox{\boxWWV}{\includegraphics[trim=0 7 0 0]{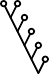}}
\def\WWV{\usebox{\boxWWV}}

\newsavebox{\boxmWWV}
\sbox{\boxmWWV}{\includegraphics[width=4mm, trim=0 5 0 0]{symbol_WWV.pdf}}
\def\mWWV{\usebox{\boxmWWV}}

\newsavebox{\boxVY}
\sbox{\boxVY}{\includegraphics[trim=0 7 0 0]{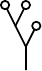}}
\def\VY{\usebox{\boxVY}}

\newsavebox{\boxIW}
\sbox{\boxIW}{\includegraphics[trim=0 5 0 0]{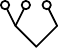}}
\def\IW{\usebox{\boxIW}}

\newsavebox{\boxVC}
\sbox{\boxVC}{\includegraphics[trim=0 5 0 0]{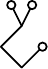}}
\def\VC{\usebox{\boxVC}}

\newsavebox{\boxCC}
\sbox{\boxCC}{\includegraphics[trim=0 3 0 0]{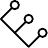}}
\def\CC{\usebox{\boxCC}}

\newsavebox{\boxBVW}
\sbox{\boxBVW}{\includegraphics[trim=0 5 0 0]{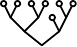}}
\def\BVW{\usebox{\boxBVW}}

\newsavebox{\boxVF}
\sbox{\boxVF}{\includegraphics[trim=0 5 0 0]{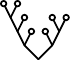}}
\def\VF{\usebox{\boxVF}}

\newsavebox{\boxVBW}
\sbox{\boxVBW}{\includegraphics[trim=0 5 0 0]{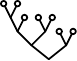}}
\def\VBW{\usebox{\boxVBW}}

\newsavebox{\boxBonB}
\sbox{\boxBonB}{\includegraphics[trim=0 5 0 0]{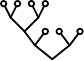}}
\def\BB{\usebox{\boxBonB}}

\newsavebox{\boxWVB}
\sbox{\boxWVB}{\includegraphics[trim=0 5 0 0]{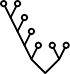}}
\def\WVB{\usebox{\boxWVB}}

\newsavebox{\boxFV}
\sbox{\boxFV}{\includegraphics[trim=0 5 0 0]{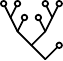}}
\def\FV{\usebox{\boxFV}}

\newsavebox{\boxBWV}
\sbox{\boxBWV}{\includegraphics[trim=0 5 0 0]{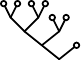}}
\def\BWV{\usebox{\boxBWV}}

\newsavebox{\boxWBV}
\sbox{\boxWBV}{\includegraphics[trim=0 5 0 0]{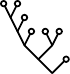}}
\def\WBV{\usebox{\boxWBV}}

\newsavebox{\boxVBVV}
\sbox{\boxVBVV}{\includegraphics[trim=0 5 0 0]{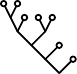}}
\def\VBVV{\usebox{\boxVBVV}}

\newsavebox{\boxBVVV}
\sbox{\boxBVVV}{\includegraphics[trim=0 5 0 0]{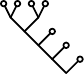}}
\def\BVVV{\usebox{\boxBVVV}}

\newsavebox{\boxWWW}
\sbox{\boxWWW}{\includegraphics[trim=0 5 0 0]{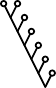}}
\def\WWW{\usebox{\boxWWW}}

\begin{document}

\title{KPZ equation with fractional derivatives of white noise}
\author{Masato Hoshino}
\date{\today}
\address{The University of Tokyo, 3-8-1 Komaba, Meguro-ku, Tokyo 153-8914, Japan}
\email{hoshino@ms.u-tokyo.ac.jp}
\renewcommand{\include}[1]{}
\renewcommand\documentclass[2][]{}
\maketitle

\begin{abstract}
In this paper, we consider the KPZ equation driven by space-time white noise replaced with
its fractional derivatives of order $\gamma>0$ in spatial variable. A well-posedness theory for the KPZ equation is established by Hairer \cite{reg} as an application of the theory of regularity structures. Our aim is to see to what extent his theory works if noises become rougher. We can expect that his theory works if and only if $\gamma<1/2$. However, we show that the renormalization like ``$(\partial_x h)^2-\infty$" is well-posed only if $\gamma<1/4$.
\end{abstract}

\documentclass{amsart}

\section{Introduction}

In this paper, we discuss the stochastic partial differential equation
\begin{align}\label{introduction:fractional kpz}
\d_th(t,x)=\d_x^2h(t,x)+(\d_xh(t,x))^2+\d_x^\g\xi(t,x)
\end{align}
for $(t,x)\in[0,\infty)\times\T$ with $\T=\R/\Z$, which is equivalent to $[0,1]$ with periodic boundary conditions, and $\g\ge0$. Here $h(t,x)$ is a continuous stochastic process and $\xi$ is a space-time white noise. $\d_x^\g=-(-\d_x^2)^\f{\g}{2}$ is the fractional derivative. If $\g=0$, (\ref{introduction:fractional kpz}) is the KPZ equation, which is proposed in \cite{KPZ} as a model of surface growth.

The equation (\ref{introduction:fractional kpz}) is ill-posed. Formally speaking, $h$ has the same regularity as the solution of the linear equation
\begin{align}\label{introduction:SHE}
\d_th(t,x)=\d_x^2h(t,x)+\d_x^\g\xi(t,x).
\end{align}
Then $h(t,\cdot)$ belongs to H\"older space $\cal{C}^\a(\T)$ with $\a<\f{1}{2}-\g$ for each fixed $t$. However this implies that the nonlinear term $(\d_xh)^2$ is the square of the distribution, which generally does not make sense.

Hairer discussed the solution of the KPZ equation in \cite{Solving KPZ} and \cite{Singular SPDE}. It is natural to replace $\xi$ by a smooth approximation $\xi_\e$, which is obtained by a convolution with a smooth mollifier, and consider the classical solution of the KPZ equation with $\xi_\e$. He showed that there exists a sequence of constants $C_\e\sim\f{1}{\e}$ such that, the sequence of solutions $h_\e$ of
\begin{align}\label{introduction:renormalized KPZ}
\d_t h_\e(t,x) = \d_x^2 h_\e(t,x) + (\d_x h_\e(t,x))^2-C_\e+\xi_\e(t,x)
\end{align}
has a unique limit $h$ in probability, which is independent of the choice of a mollifier.

Our goal is to make the noise rougher and see to what extent this theory works. Because of the ``local subcriticality" (Assumption 8.3 of \cite{reg}), we can expect that similar results hold if $\g<\f{1}{2}$. If we write $h^\delta(t,x)=\delta^{-\f{1}{2}+\g}h(\delta^2t,\delta x)$ and $\xi^\delta(t,x)=\delta^{\f{3}{2}}\xi(\delta^2t,\delta x)$ for $\delta>0$, then $\xi^\delta$ is equal to $\xi$ in distribution and $h^\delta$ satisfies
\[\d_th^\delta(t,x)=\d_x^2h^\delta(t,x)+\delta^{\f{1}{2}-\g}(\d_xh^\delta(t,x))^2+\d_x^\g\xi^\delta(t,x).\]
As $\delta\to0$, we can see that the nonlinear term vanishes. Formally speaking, this means that $h$ behaves like the solution of (\ref{introduction:SHE}) at small scales. His theory implies that it is possible to devise a suitable renormalization in this case. However, we prove that the renormalization like (\ref{introduction:renormalized KPZ}) is possible only if $\g<\f{1}{4}$ in this paper. In the case $\g\ge\f{1}{4}$, see Subsection \ref{Subsection:case of 1/4}.

\begin{theo}\label{Theorem:main result}
Let $\rho$ be a function on $\R^2$ which is smooth, compactly supported, symmetric in $x$, nonnegative, and satisfies $\int_{\R^2}\rho(t,x)dtdx=1$. Set $\rho_\e(t,x)=\e^{-3}\rho(\e^{-2}t,\e^{-1}x)$ for $\e>0$, and $\xi_\e=\xi*\rho_\e$ (space-time convolution). Let $0\le\g<\f{1}{4}$ and $0<\eta<\f{1}{2}-\g$. Then there exists a sequence of constants $C_\e\sim_{\g,\rho}\e^{-1-2\g}$ such that, for every initial condition $h_0\in\cal{C}^\a(\T)$ the sequence of solutions $h_\e$ of
\[\d_th_\e(t,x)=\d_x^2h_\e(t,x)+(\d_xh_\e(t,x))^2-C_\e+\d_x^\g\xi_\e(t,x)\]
converges to a unique stochastic process $h$, which is independent of the choice of $\rho$. Precisely, $h_\e(t,\cdot)$ exists until the survival time $T_\e\in(0,\infty]$ which satisfies $\liminf_{\e\downarrow0}T_\e>0$, and $h_\e$ convergences to $h$ in probability in the uniform norm on $[0,T]\times\T$ and $\eta$-H\"older norm on all compact sets in $(0,T]\times\T$ for every $T<\liminf_\e T_\e$.
\end{theo}

This theorem is obtained from Theorem \ref{Theorem:solution map}, Proposition \ref{Proposition:renormalized KPZ for 1/6} and Theorem \ref{Theorem:convergence of model for 1/6} for $0\le\g<\f{1}{6}$, and from Theorem \ref{Theorem:solution map}, Proposition \ref{Proposition:renormalized KPZ for 1/4} and Theorem \ref{Theorem:convergence of model for 1/4} for $\f{1}{6}\le\g<\f{1}{4}$. The estimate of $C_\e$ is in Propositions \ref{Prop:const1} and \ref{Prop:const2}.

We note that above result is related to \cite{HHLNT}, where the equation
\[\d_tu(t,x)=\d_x^2u(t,x)+u\d_t\d_xW(t,x)\]
is studied. Here $W$ is a standard Brownian motion in $t$ and a fractional Brownian motion with Hurst parameter $H\in(\f{1}{4},\f{1}{2})$ in $x$. The relation between both equations is $H=\f{1}{2}-\gamma$, so that the same boundary $\gamma=H=\f{1}{4}$ appears. Although both equations are same after performing the Cole-Hopf transformation $u=e^h$, their results do not imply Theorem \ref {Theorem:main result} since It\^o calculus is only stable under spatial regularizations.

The organization of this paper is as follows. In Section \ref{Section:notation}, we introduce some notations and fractional calculus. In Section \ref{Section:reg}, we briefly recall the theory of regularity structures and prepare some tools for the proof of Theorem \ref{Theorem:main result}. We discuss the renormalization of models in Section \ref{Section:renormalization}. Details of the proof are given in Section \ref{Section:estimate1} ($0\le\g<\f{1}{10}$), Section \ref{Section:estimate2} ($\f{1}{10}\le\g<\f{1}{6}$), Section \ref{Section:estimate3} ($\f{1}{6}\le\g<\f{3}{14}$), and Section \ref{Section:estimate4} ($\f{3}{14}\le\g<\f{1}{4}$). 

\end{document}